\documentclass[12pt,leqno,draft]{article}

\newtheorem{theorem}{Theorem}
\newtheorem{lemma}[theorem]{Lemma}
\newtheorem{proposition}[theorem]{Proposition}
\newtheorem{definition}[theorem]{Definition}
\newtheorem{corollary}[theorem]{Corollary}

\newcommand{\begintheorem}{\addtocounter{equation}{1}\begin{theorem}}
\newcommand{\beginlemma}{\addtocounter{equation}{1}\begin{lemma}}
\newcommand{\beginproposition}{\addtocounter{equation}{1}\begin{proposition}}
\newcommand{\begindefinition}{\addtocounter{equation}{1}\begin{definition}}
\newcommand{\begincorollary}{\addtocounter{equation}{1}\begin{corollary}}

\begin{document}

\title{Notes on normed algebras, 2}

\author{Stephen William Semmes	\\
	Rice University		\\
	Houston, Texas}

\date{}

\maketitle

	Let $\mathcal{A}$ be a finite-dimensional commutative algebra
over the complex numbers, with identity element $e$.  Thus
$\mathcal{A}$ is a finite-dimensional complex vector space equipped
with an additional binary operation of multiplication which satisfies
the usual rules of associativity, commutativity, and distributivity,
and $e$ is a nonzero element of $\mathcal{A}$ such that $e \, x = x$
for all $x \in \mathcal{A}$.  As a basic class of examples, one can
take $X$ to be a finite set, and $\mathcal{A}$ the algebra of
complex-valued functions on $X$, with respect to pointwise
multiplication.

	As another class of examples, suppose that $A$ is a
commutative semigroup with identity.  In other words, $A$ is a set
equipped with an identity element $0$ and a binary operation $+$ which
satisfies the usual associativity and commutativity rules, and with $x
+ 0 = x$ for all $x \in A$.  Let $\mathcal{A}$ be the vector space of
complex-valued functions on $A$, and define the convolution of two
such functions $f_1$, $f_2$ by
\begin{equation}
	(f_1 * f_2)(z) = \sum_{x + y = z} f_1(x) \, f_2(y).
\end{equation}
More precisely, the sum is taken over all $x, y \in A$ such that $x +
y = z$.  One can check that this operation of convolution satisfies
the commutative and associative laws, and of course it is linear in
$f_1$, $f_2$.

	For each $x \in A$, let $\delta_x$ denote the function on $A$
which is equal to $1$ at $x$ and to $0$ at all other elements of $A$.
The semigroup operation on $A$ corresponds exactly to convolution of
these functions.  The function $\delta_0$ serves as an identity
element for the operation of convolution.  Thus the vector space of
complex-valued functions on $A$ becomes a commutative algebra with
respect to convolution.

	Let $\mathcal{A}$ be a finite-dimensional commutative algebra
over the complex numbers.  By an ideal in $\mathcal{A}$ we mean a
linear subspace $\mathcal{I}$ of $\mathcal{A}$ such that if $x \in
\mathcal{I}$ and $a \in \mathcal{A}$, then $a \, x \in \mathcal{I}$ as
well.  We say that an ideal $\mathcal{I}$ in $\mathcal{A}$ is proper
if it is not all of $\mathcal{A}$.  As a special case, if $x \in
\mathcal{A}$, then we get an ideal $\mathcal{I}(x)$ consisting of $a
\, x$ as $a$ runs through all elements of $\mathcal{A}$.  This ideal
is proper if and only if $x$ does not have an inverse in
$\mathcal{A}$.

	Suppose that $\mathcal{I}$ is a proper ideal in the
commutative finite-dimensional algebra $\mathcal{A}$ over the complex
numbers.  Thus $\mathcal{I}$ is a linear subspace of $\mathcal{A}$,
and we can form the quotient space $\mathcal{A} / \mathcal{I}$
initially as a vector space of positive dimension equal to the
dimension of $\mathcal{A}$ minus the dimension of $\mathcal{I}$.
There is an associated quotient mapping from $\mathcal{A}$ onto
$\mathcal{A} / \mathcal{I}$, which sends an element of $\mathcal{A}$
to the corresponding element of the quotient space.  The kernel of
this mapping is exactly equal to $\mathcal{I}$.

	By the usual arguments, the operation of multiplication on
$\mathcal{A}$ leads to a similar operation on the quotient
$\mathcal{A} / \mathcal{I}$, because $\mathcal{I}$ is an ideal.  This
operation on the quotient $\mathcal{A} / \mathcal{I}$ is commutative,
because the initial operation on $\mathcal{A}$ is commutative.  The
image of the identity element $e$ in the quotient is an identity
element in the quotient, and it is nonzero because the ideal
$\mathcal{I}$ is proper, and therefore does not contain $e$.
In short the quotient $\mathcal{A} / \mathcal{I}$ is itself a
finite-dimensional commutative algebra.

	A maximal ideal $\mathcal{I}$ in $\mathcal{A}$ is a proper
ideal with the property that any ideal in $\mathcal{A}$ which contains
$\mathcal{I}$ is either equal to $\mathcal{I}$ or to $\mathcal{A}$.
Every proper ideal in $\mathcal{A}$ is contained in a maximal ideal,
namely, in a proper ideal of maximal dimension.  The only proper ideal
in $\mathcal{A}$ is the ideal containing only $0$ if and only if
$\mathcal{A}$ is a field, which is to say that every nonzero element
of $\mathcal{A}$ has an inverse in $\mathcal{A}$.  A finite-dimensional
algebra over the complex numbers which is a field is equal to the
span of its identity element, which is to say that it has dimension
$1$ and is isomorphic to the complex numbers.  This follows from the
fundamental theorem of algebra.

	Suppose that $\mathcal{A}$ is a finite-dimensional commutative
algebra over the complex numbers with identity element $e$, and that
$\mathcal{I}$ is a proper ideal in $\mathcal{A}$.  Thus we get a
nonzero homomorphism from $\mathcal{A}$ onto the quotient $\mathcal{A}
/ \mathcal{I}$.  If $\mathcal{I}$ is a maximal ideal, then
$\mathcal{A} / \mathcal{I}$ is one-dimensional and isomorphic to the
complex numbers.  Equivalently, we get a homomorphism from
$\mathcal{A}$ onto the complex numbers whose kernel is equal to
$\mathcal{I}$ and which sends the identity element $e$ in
$\mathcal{A}$ to the complex number $1$.  Conversely, a homomorphism
from $\mathcal{A}$ onto the complex numbers has kernel equal to a
maximal ideal in $\mathcal{A}$ and sends $e$ to $1$.

	Suppose that $\mathcal{A}$ is a finite-dimensional commutative
algebra over the complex numbers with identity element $e$, and that
$\phi$ is a homomorphism from $\mathcal{A}$ onto the complex numbers,
which then takes $e$ to $1$ automatically.  If $x$ is an invertible
element of $\mathcal{A}$, then $\phi(x)$ is an invertible complex
number, which means that $\phi(x) \ne 0$.  Thus $\phi(x) = 0$ implies
that $x$ does not have an inverse in $\mathcal{A}$.  If we start with
$x \in \mathcal{A}$ which is not invertible, then the ideal
$\mathcal{I}(x)$ generated by $x$ is proper, and is therefore
contained in a maximal ideal.  It follows that there is a homomorphism
$\phi$ from $A$ onto the complex numbers such that $\phi(x) = 0$.

	Recall that the spectrum of an element $x$ of $\mathcal{A}$ is
defined to be the set of complex numbers $\lambda$ such that $\lambda
\, e - x$ does not have an inverse in $\mathcal{A}$.  The preceding
remarks imply that the spectrum is the same as the set of complex
numbers which arise as $\phi(x)$ for some homomorphism $\phi$ from
$\mathcal{A}$ onto the complex numbers.  We also know that the
spectrum of $x$ is a nonempty finite set of complex numbers.

	Now suppose that $\mathcal{A}$ is a finite-dimensional
commutative algebra over the complex numbers equipped with a norm
$\|\cdot \|$ such that $(\mathcal{A}, \|\cdot \|)$ is a normed
algebra.  Thus $\|\cdot \|$ is a nonnegative real-valued function on
$\mathcal{A}$ such that $\|x\| = 0$ if and only if $x = 0$, $\|\alpha
\, x\| = |\alpha| \, \|x\|$ for all complex numbers $\alpha$ and $x
\in \mathcal{A}$, $\|e\| = 1$, and the norm of a sum or product of two
elements of $\mathcal{A}$ is less than or equal to the corresponding
product of the norms.  If $\lambda$ is a complex number and $x$ is an
element of $\mathcal{A}$ such that $\|x\| < |\lambda|$, then $\lambda
\, e - x$ is invertible in $\mathcal{A}$, because one can sum the
series $\sum_{j=0}^\infty \lambda^{-j} \, x^j$.  In other words,
$|\lambda| \le \|x\|$ for all $\lambda$ in the spectrum of $x$.  If
$\phi$ is a homomorphism from $\mathcal{A}$ onto the complex numbers,
then $|\phi(x)| \le \|x\|$ for all $x \in \mathcal{A}$, since
$\phi(x)$ is in the spectrum of $\mathcal{A}$.

	Let $X$ be a finite nonemtpy set, and let $\mathcal{A}$ be the
algebra of complex-valued functions on $X$, with respect to pointwise
addition and multiplication.  If $\mathcal{I}$ is an ideal in
$\mathcal{A}$, then one can show that there is a subset $E$ of
$X$ such that $\mathcal{I}$ consists exactly of the functions
$f$ on $X$ which satisfy $f(p) = 0$ when $p \in E$.  Moreover,
$\mathcal{I}$ is a proper ideal if and only if there is a
point $p \in X$ such that $\mathcal{I}$ consists exactly of the
functions $f$ on $X$ such that $f(p) = 0$.  Thus the homomorphisms
from $\mathcal{A}$ onto the complex numbers are exactly of the
form $f \mapsto f(p)$ for $p \in X$.

	There is a natural norm in this setting, defined by saying
that the norm of a function $f$ on $X$ is equal to the maximum of
$|f(p)|$ over $p \in X$.  The spectrum of a function $f$ on $X$, as an
element of the algebra of functions on $X$, is equal to the set of
values of $f$.  Thus the norm of $f$ is equal to the maximum of the
absolute values of the complex numbers in the spectrum of $f$.

	Let $A$ be a finite commutative semigroup, and let
$\mathcal{A}$ be the algebra of complex-valued functions on $A$ with
respect to convolution as described earlier.  Suppose that $\phi$ is a
homomorphism from $\mathcal{A}$ onto the complex numbers.  This leads
to a mapping $\Phi$ from $A$ to the complex numbers, defined by
$\Phi(a) = \phi(\delta_a)$ for $a \in A$, where $\delta_a$ is the
function on $A$ which is equal to $1$ at $a$ and to $0$ at other
elements of $A$.  When $a = 0$, this is the identity element of
$\mathcal{A}$.  The requirement that $\phi$ map $\mathcal{A}$ onto the
complex numbers is equivalent to $\phi$ mapping $\delta_0$ to $1$,
which is the same as $\Phi(0) = 1$.

	Because $\phi : \mathcal{A} \to {\bf C}$ is a homomorphism,
$\Phi$ should be a homomorphism from $A$ into ${\bf C}$ as a semigroup
with respect to multiplication.  In other words we should have $\Phi(a
+ b) = \Phi(a) \, \Phi(b)$ for all $a, b \in A$.  Under the assumption
that $\phi$ is linear, this condition on $\Phi$ is equivalent to $\phi
(f_1 * f_2) = \phi(f_1) \, \phi(f_2)$ for all $f_1, f_2 \in
\mathcal{A}$.  This follows from the fact that $\mathcal{A}$ is
spanned by the functions $\delta_a$, $a \in A$.

	Let us define a norm on $\mathcal{A}$ in this situation by
setting $\|f\|$ equal to $\sum_{a \in A} |f(a)|$ for all
complex-valued functions on $A$.  One can check that 
\begin{equation}
	\|f_1 * f_2\| \le \|f_1\| \, \|f_2\|
\end{equation}
for all complex-valued functions $f_1$, $f_2$ on $A$ in this
situation.  Also, $\|\delta_a\| = 1$ for all $a \in A$, and in
particular the norm of the identity element of $\mathcal{A}$ is equal
to $1$.

	Suppose that $\phi$ is a homomorphism from $\mathcal{A}$ onto
the complex numbers, and let $\Phi$ be the corresponding homomorphism
from $A$ into the multiplicative semigroup of complex numbers.  The
image of $A$ under $\Phi$ is a finite subsemigroup of the
multiplicative semigroup of complex numbers, and therefore $|\Phi(a)|
\le 1$ for all $a \in A$.  Indeed, if we choose $a \in A$ so that
$|\Phi(a)|$ is maximal, then $|\Phi(a + a)|$ is equal to
$|\Phi(a)|^2$, and it is less than or equal to $|\Phi(a)|$ by
maximality.  Thus $|\Phi(a)| \le 1$, as desired.  Of course $\Phi(0) =
1$, so that the maximum of $|\Phi|$ on $A$ is equal to $1$.  If $f$ is
a complex-valued function on $A$, then we can use linearity of $\phi$
to write $\phi(f)$ explicitly as $\sum_{a \in A} \Phi(a) \, f(a)$.  We
can apply absolute values to this sum and use the triangle inequality
to obtain directly that $|\phi(f)| \le \|f\|$.


\begin{thebibliography}{1}

\bibitem {1} W.~Rudin, {\it Functional Analysis}, second edition,
McGraw-Hill, 1991.




\end{thebibliography}
\end{document}